# UPCROSSING INEQUALITIES FOR STATIONARY SEQUENCES AND APPLICATIONS

By Michael Hochman

*Hebrew University of Jerusalem*

For arrays $(S_{i,j})_{1 \leq i \leq j}$ of random variables that are stationary in an appropriate sense, we show that the fluctuations of the process $(S_{1,n})_{n=1}^\infty$ can be bounded in terms of a measure of the "mean subadditivity" of the process $(S_{i,j})_{1 \leq i \leq j}$. We derive universal upcrossing inequalities with exponential decay for Kingman's subadditive ergodic theorem, the Shannon–MacMillan–Breiman theorem and for the convergence of the Kolmogorov complexity of a stationary sample.

**1. Introduction.** Let us say that a sequence $(X_n)_{n=1}^\infty$ of real numbers has $k$ crossings (or upcrossings) of an interval $[s,t]$ if there are indices

$$1 \leq i_1 < j_1 < i_2 < j_2 < \cdots < i_k < j_k$$

such that $X_{i_m} < s$ and $X_{j_m} > t$. Allowing $X_n$ to be random, it easily follows that $\lim X_n$ exists a.s. if and only if, for every interval of positive length, the probability of infinitely many crossings of the interval is 0.

There are a number of classical limit theorems in probability that can be formulated and proven in this way, the best known of which is Doob's upcrossing inequality for $L^1$ martingales [6]: if $(S_n)_{n=1}^\infty$ is an $L^1$ martingale, then, for $s < t$,

$$\mathbb{P}((S_n)_{n=1}^\infty \text{ has } k \text{ upcrossings of } [s,t]) \leq \frac{\sup_n \|S_n\|_1}{k(t-s)}$$

(see also Dubins [7]). A similar inequality was proven by Bishop for the time averages $S_n = \frac{1}{n}\sum_{i=1}^n X_i$ of an $L^1$ stationary process $(X_n)_{n=1}^\infty$ [1, 2]. Assuming nonnegativity of the process instead of integrability, Ivanov [8] proved the following beautiful result: for every $s < t$,

(1.1) $$\mathbb{P}((S_n)_{n=1}^\infty \text{ has } k \text{ upcrossings of } [s,t]) \leq \left(\frac{s}{t}\right)^k$$

(see [4]; for related results, see Jones et al. [9], and Kalikow and Weiss [10]). A remarkable aspect of these inequalities is that they hold universally: except for trivial normalization, they do not depend on the process in question. Neither martingales nor ergodic averages admit universal rates of convergence and it is all the more surprising that such general bounds for the fluctuations exist.

In this paper, we establish a general upcrossing inequality for certain sequences associated with stationary processes, in terms of a certain measure of "mean subadditivity" of the process. We consider arrays $(S_{i,j})_{1 \leq i \leq j}$ of random variables, although we will usually identify the ordered pair $(i,j)$ that indexes $S_{i,j}$ with the integer interval $[i;j] = [i,j] \cap \mathbb{N}$, and if $U = [i;j]$ is such an interval, we write $S_U$ for $S_{i,j}$. We assume the process to be stationary with respect to translation of the indexing intervals, that is, that for every $m \in \mathbb{N}$, one has

$$(S_{[i;j]})_{1 \leq i \leq j} = (S_{[i+m,j+m]})_{1 \leq i \leq j} \qquad \text{in distribution.}$$

One very general way to obtain such arrays is by applying a function to samples of stationary processes: if $(X_n)_{n=1}^\infty$ is stationary and $g$ is any function defined on finite sequences, then $S_{i,j} = g(X_i, X_{i+1}, \ldots, X_j)$ satisfies these assumptions.

Let $I \subseteq \mathbb{N}$ be an interval and $\delta > 0$. We say that a collection $I_1, \ldots, I_r$ of intervals $\delta$-*fills* $I$ if all of the intervals are contained in $I$ and $|I \setminus \bigcup I_i| < \delta |I|$, where $|\cdot|$ denotes cardinality.

THEOREM 1.1. *Suppose that $(S_{i,j})_{1 \leq i \leq j}$ is stationary in the above sense. Let $s < t$ and $0 < \delta < \frac{1}{4}$. Then, for every $k$,*

(1.2) $$\mathbb{P}\begin{pmatrix} (S_{1,n})_{n=1}^\infty \text{ has } k \\ \text{upcrossing of } [s,t] \end{pmatrix}$$
$$\leq c \cdot \rho^k + \mathbb{P}\begin{pmatrix} \text{there exists } n > k \text{ such that } S_{1,n} > t \\ \text{and } [1;n] \text{ can be } \delta\text{-filled by disjoint} \\ \text{intervals } V_1, \ldots, V_r \text{ satisfying } S_{V_j} < s \end{pmatrix}.$$

*The constants $c$ and $0 < \rho < 1$ depend only on $\delta$ (but not on the process or on $s, t$).*

In applications, one optimizes over $\delta$ to get a bound for the left-hand side which is often independent of the process. Theorem 1.1 is effective and the constants may be computed explicitly, although they are surely not optimal.

This inequality cannot be reversed and in Section 2, we give a simple example in which the left-hand side decays exponentially and uniformly for a certain class of processes, but the right-hand side can decay arbitrarily slowly.

Theorem 1.1 can be generalized in several ways. It remains valid when one starts with a "two-sided" stationary array $(S_{i,j})_{-\infty<i\leq j<\infty}$ and sets $S_n = S_{-n,n}$, and there is a version for $S_n = S_{V_n}$, where $V_i$ is an arbitrary increasing sequence with $0 \in V_1$. It can also be extended to the multidimensional setting, where the process is indexed by cubes instead of segments. These versions require minor modifications of the proof we give.

Our first application is to Kingman's subadditive ergodic theorem. Let $(X_{m,n})_{1 \leq m \leq n}$ be stationary, in the sense described above, and subadditive, that is, $X_{k,n} \leq X_{k,m} + X_{m+1,n}$ whenever $k \leq m < n$. Kingman's theorem states that under some integrability conditions, if $X_{m,n}$ is stationary and subadditive, then $\frac{1}{n} X_{1,n}$ converges almost surely. As examples of this situation, consider $X_{m,n} = \sum_{i=m}^{n} Y_i$, where $Y_i$ is a stationary process or $X_{m,n} = \|Z_m, \ldots, Z_n\|$, where $Z_i$ is a stationary sequence of operators.

Polynomial decay for upcrossings in Kingman's theorem for integrable processes was proven by Krawczak [11]. We establish an exponential version of this, with integrability replaced by a boundedness condition:

THEOREM 1.2. *Let $(X_{m,n})_{1 \leq m \leq n}$ be stationary and subadditive. Suppose that $X_{1,1} \leq M$ a.e. for a constant $M$. Then, for every $s < t$,*

$$\mathbb{P}((S_n)_{n=1}^{\infty} \text{ has } k \text{ upcrossings of } [s,t]) < c \cdot \rho^k$$

*for constants $0 < \rho < 1$ and $c$ that depend only on $s, t$ and $M$.*

We next turn to the convergence of Shannon information of samples. Given a finite-valued process $(X_i)_{i=1}^{\infty}$, the information of the sample $X_1, \ldots, X_n$ is the random variable

$$I(X_1, \ldots, X_n) = -\log \mathbb{P}(X_1, \ldots, X_n).$$

Here, for a fixed sequence $\xi = \xi_1, \ldots, \xi_n$, we write $\mathbb{P}(\xi)$ for the probability of observing this sample, so $\mathbb{P}(X_1, \ldots, X_n)$ is the probability of observing the sample that was actually observed.

The Shannon–MacMillan–Breiman theorem [13] is one of the fundamental theorems in information theory, asserting that $\frac{1}{n} I(X_1, \ldots, X_n)$ converges almost surely. See, for instance, [12].

THEOREM 1.3. *Suppose that $(X_n)_{n=1}^{\infty}$ is a $0,1$-valued stationary process. Let $S_n = \frac{1}{n} I(X_1, \ldots, X_n)$. Then, for every $s < t$,*

$$\mathbb{P}((S_n)_{n=1}^{\infty} \text{ has } k \text{ upcrossings of } [s,t]) < c \cdot \rho^k$$

*and the constants $0 < \rho < 1$ and $c$ depend on $s, t$.*

The same result holds for processes with a finite number $r$ of symbols, but the constants will then also depend on $r$. We note that there is no universal rate of convergence for this limit and Theorem 1.3 seems to be the first effective proof of its convergence. It is interesting that the classical proofs of the convergence of $\frac{1}{n}I(X_1,\ldots,X_n)$ rely on the ergodic theorem and the Martingale theorem, but despite upcrossing inequalities being known for both of these, it is unclear how to combine them to deduce Theorem 1.3.

Bishop's upcrossing proof of the ergodic theorem was part of a larger program of Bishop's to "constructivize" mathematical analysis. As Theorem 1.1 is effective, it can probably be adapted to this setting, leading to constructive proofs of the Kingman and Shannon–MacMillan–Breiman theorems.

Another closely related result involves algorithmic complexity rather than Shannon information. For a finite string $x$ of 0's and 1's, let $\kappa(x)$ denote the Kolmogorov complexity of $x$, sometimes referred to as the *minimal description length* of $x$ (see Section 4). Although $\kappa(\cdot)$ is not formally computable, it has been extensively studied as a nonstatistical measure of complexity and is closely related to Shannon's theory of information; see, for example, [5]. The following result is therefore an analog of Theorem 1.3; the existence of the limit was shown by Brudno [3].

THEOREM 1.4. *Suppose that $(X_n)_{n=1}^\infty$ is a stationary $0,1$-valued process. Then, writing $S_n = -\frac{1}{n}\kappa(X_1,\ldots,X_n)$, for every $s < t$, we have*

$$\mathbb{P}((S_n)_{n=1}^\infty \text{ has } k \text{ upcrossings of } [s,t]) < c \cdot \rho^k$$

*for constants $c$ and $0 < \rho < 1$ that depend only on $s,t$.*

The rest of the paper is organized as follows. First, we use Theorem 1.1 to derive the applications. In Section 2, we derive the result on Kingman's theorem, and show that the inequality in Theorem 1.1 cannot be reversed. In Section 3, we prove the upcrossing inequality for the Shannon–McMillan–Breiman theorem and in Section 4, we derive the inequality for Kolmogorov complexity. In Section 5, we reduce Theorem 1.1 to a combinatorial lemma whose proof is given in Section 6.

**2. The subadditive ergodic theorem and an example.** In this section, we discuss the relation between Theorem 1.1 and ergodic theorems.

PROOF OF THEOREM 1.2. Suppose that $V_1,\ldots,V_r$ are disjoint subintervals of $[1;n]$ such that $\frac{1}{|V_i|}X_{V_i} < s$. Let $U = [0,1] \setminus \bigcup_{i=1}^r V_i$. By subadditivity,

$$\frac{1}{n}X_{1,n} \leq \sum \frac{|V_i|}{n}\frac{1}{|V_i|}X_{V_i} + \sum_{j \in U} X_{j,j} \leq \frac{\sum |V_i|}{n}s + \frac{|U|}{n}M$$

(in the last inequality, we used stationarity to get $X_{j,j} \leq M$). Thus, if $[1;n]$ is $\delta$-filled by the $V_i$, then $|U| \leq \delta n$, so

$$\frac{1}{n} X_{1,n} \leq (1-\delta)s + \delta M$$

and if $\delta = \frac{t-s}{M-s}$, then it impossible that $\frac{1}{n} X_{1,n} > t$. Hence, for this $\delta$, when we apply Theorem 1.1 with the sequence $S_{i,j} = \frac{1}{j-i+1} X_{i,j}$, the event on the right-hand side of inequality (1.2) is empty. The theorem follows. □

Next, we show that decay of the right-hand side in Theorem 1.1 is not necessary in order to get fast decay on the left. Let $(X_n)_{n=1}^{\infty}$ be a stationary process with values in $[-1, 1]$ and let

$$S_{i,j} = \frac{1}{\lfloor \sqrt{j-i} \rfloor} \sum_{k=i}^{i+\lfloor \sqrt{j-i} \rfloor - 1} X_k,$$

so $S_{1,n}$ is obtained by repeating elements from the sequence of ergodic averages $(\frac{1}{n}(X_1 + \cdots + X_n))_{n=1}^{\infty}$. The latter sequence obeys an upcrossing inequality which does not depend on $X_n$; thus, the former sequence also does. However, the following proposition shows that one cannot obtain this from Theorem 1.1.

PROPOSITION 2.1. *Let $\delta > 0$, and*

$$p_k = \mathbb{P} \begin{pmatrix} \text{there exists } n > k \text{ such that } S_{1,n} > \frac{1}{2} \\ \text{and } [1,n] \text{ can be } \delta\text{-filled by disjoint} \\ \text{intervals } V_1, \ldots, V_r \text{ satisfying } S_{V_j} < -\frac{1}{2} \end{pmatrix}.$$

*Then, for arbitrarily large $k$, there are processes for which $p_k \geq 1/6$ and, in particular, the convergence $p_k \to 0$ is not uniform as the process $(X_n)$ is varied.*

PROOF. Fix $\delta > 0$. Let $n > 1/\delta$ and let $(X_i)_{i=1}^{\infty}$ be the process whose unique sample path, up to translation, is the sequence with period $2n$ in which blocks of 1's and $-1$'s of length $n$ alternate. Set $k = n^2$. It is easily verified that (a) the probability that the first $2n/3$ symbols of a sample are 1's is $1/6$ and $S_{1,k}$, in this case, is $\geq 1/2$; and (b) if $X_1 = 1$, then, taking $j$ to be the first index with $X_j = -1$, we have $S_{j,k} < -1/2$ and $[1;k]$ is $\delta$-covered by $[j;k]$. Thus, for each square $k$, there are processes for which $p_k \geq 1/6$. □

**3. The Shannon–McMillan–Breiman theorem.** In this section, we prove Theorem 1.3. Fix $s < t$ and a parameter $\delta > 0$. Set

$$S_{i,j} = \frac{1}{j-i+1} I(X_i, \ldots, X_j).$$

In order to apply Theorem 1.1, for each $n \in \mathbb{N}$, we wish to bound the probability of the event
$$B_n = \left\{ \begin{array}{c} S_{1,n} > t \text{ and } [1,n] \text{ can be } \delta\text{-filled by disjoint} \\ \text{intervals } V_1, \ldots, V_r \text{ satisfying } S_{V_j} < s \end{array} \right\}.$$

Since
$$S_{i,j} > t \iff \mathbb{P}(X_i, \ldots, X_j) < 2^{-t(j-i+1)},$$
$$S_{i,j} < s \iff \mathbb{P}(X_i, \ldots, X_j) > 2^{-s(j-i+1)},$$

we have the trivial bound

$$(3.1) \quad \mathbb{P}(B_n) \leq 2^{-tn} \# \left( \begin{array}{c} \text{words } w \in \{0,1\}^n \text{ which can be } \delta\text{-filled with} \\ \text{disjoint words } v \text{ satisfying } \mathbb{P}(v) > 2^{-s\ell(v)} \end{array} \right),$$

where $\ell(v)$ denotes the length of $v$.

To estimate the right-hand side of (3.1), we note that each word $w$ that is counted on the right can be constructed as follows:

(1) choose a subset $I \subseteq [1;n]$ of size $\leq \delta n$;
(2) choose the symbol $w(i)$ for each $i \in I$;
(3) for each maximal interval $J \subseteq [1;n] \setminus I$, choose a word $v = w|_J$ with $\mathbb{P}(v) > 2^{-s\ell(v)}$.

To bound the number of words produced in (1)–(3), we bound the number of choices at each step. In step (1), we have $\leq 2^{nh(\delta)+o(\log n)}$ choices, where $h(x) = -x \log x - (1-x) \log(1-x)$ (this is a standard consequence of Stirling's formula). In step (2), we have at most $2^{\delta n}$ choices. Finally, in step (3), let $J_1, \ldots, J_r$ be the maximal intervals in $[1;n] \setminus I$. The number of distinct words $v$ of length $|J_i|$ and satisfying $\mathbb{P}(v) > 2^{-s|J_i|}$ is clearly bounded by $2^{s|J_i|}$, so the number of ways to choose such words with lengths $|J_1|, \ldots, |J_r|$ is at most
$$\prod_{i=1}^{r} 2^{s|J_i|} = 2^{\sum_{i=1}^{r} s|J_i|} \leq 2^{sn}.$$

It follows that the number of words counted on the right-hand side of equation (3.1) is
$$\leq 2^{(s+h(\delta)+\delta)n+o(\log n)},$$

so
$$\mathbb{P}(B_n) \leq 2^{-(t-s-h(\delta)-\delta+o(1))n}.$$

Hence, if $\delta$ is small enough in a manner depending on $s, t$, then this bound is summable and gives
$$\mathbb{P}\left( \bigcup_{n > k} B_n \right) < c \cdot \rho^k$$

for constants $c$ and $0 < \rho < 1$ depending only on $s, t$. Since $\bigcup_{n>k} B_n$ is the event on the right-hand side of inequality (1.1) in Theorem 1.1, this completes the proof of Theorem 1.3.

**4. Kolmogorov complexity.** The Kolmogorov complexity $\kappa(x)$ of a 0,1-valued string $x$ is defined as follows. Fix a universal Turing machine $U$ and let $x^*$ be a string of minimal length such that, when $U$ is run on input $x^*$, the output is $x$. Then, $\kappa(x)$ is the length of $x^*$. Although $\kappa(\cdot)$ depends on the universal machine $U$, changing $U$ only changes $\kappa(\cdot)$ by an additive constant.

In this section, we prove Theorem 1.4. Note that $S_{m,n} = \frac{1}{n}\kappa(x_m, \ldots, x_n)$ is not strictly subadditive, so Kingman's theorem does not apply; nonetheless, we can exploit the fact that it is "almost subadditive."

PROPOSITION 4.1. *Let $(X_n)_{n=1}^\infty$ be a 0,1-valued stationary process and set $S_n = \frac{1}{n}\kappa(X_1, \ldots, X_n)$. Then, for every $s < t$, there is a $\delta > 0$ and an $N$ such that the event*

$$B_k = \left\{ \begin{array}{c} \text{there exists } n > k \text{ such that } S_{1,n} > t \\ \text{and } [1, n] \text{ can be } \delta\text{-filled by disjoint} \\ \text{intervals } V_1, \ldots, V_r \text{ satisfying } S_{V_j} < s \end{array} \right\}$$

*is empty for all $k \geq N$.*

PROOF. We reason as in the previous section. Suppose that $x \in \{0,1\}^n$ and $[1;n]$ can be $\delta$-filled by a disjoint collection of intervals $\{V_1, \ldots, V_m\}$ with $\frac{1}{|V_k|}\kappa(x|_{V_k}) < s$. We can encode $x$ by describing the choices of the three-step process outlined in the previous section. In step (3), we encode the pattern $x|_{V_m}$ by writing down the algorithm that produces it. Thus, the estimate from the previous section shows that the number of bits required is

$$h(\delta)n + n + \sum \kappa(x|_{V_i}) \leq (s + h(\delta) + \delta)n + o(\log n).$$

We only require a constant-length program to extract the string $x$ from this encoding. We also require an overhead of $O((h(\delta) + \delta)n)$ to encode this information in a self-punctuating way. Thus,

$$\kappa(x) \leq o(\log n) + (s + C'(h(\delta) + \delta))n$$

and if $\delta$ is small enough (in a manner depending on $s, t$), then this implies $\kappa(x) < t$, once $n$ is large enough (how large $n$ must be depends on $C, C'$, which, in turn, depend on the Turing machine we are using, but is independent of $s, t, \delta$ and the process). The proposition follows. □

Theorem 1.4 now follows from Theorem 1.1.

**5. Reduction of Theorem 1.1 to a covering lemma.** The purpose of this section is to reduce the proof of Theorem 1.1 to a combinatorial statement about intervals, related to the effective Vitali covering lemma of Kalikow and Weiss [10]. This lemma is stated below in Lemma 5.1, but its proof is rather technical and we defer it to Section 6.

Let $(S_{i,j})_{1 \leq i \leq j}$ be stationary, in the sense discussed in the Introduction. Fix $s < t$, a parameter $\delta > 0$ and an integer $k$. For $i \in \mathbb{N}$, we define the events

$$A_i = \{(S_{i,i+n})_{n=1}^\infty \text{ has } k \text{ upcrossings of } [s,t]\},$$

$$B_i = \left\{ \begin{array}{l} \text{There is an } n > k \text{ such that } S_{i,i+n} > t, \text{ and } [i, i+n] \text{ can be} \\ \delta\text{-filled by disjoint intervals } V_1, \ldots, V_r \text{ satisfying } S_{V_j} < s \end{array} \right\}.$$

By stationarity, $\mathbb{P}(A_i) = \mathbb{P}(A_j)$ and $\mathbb{P}(B_i) = \mathbb{P}(B_j)$ for all $i, j$. We abbreviate $A = A_1, B = B_1$. Theorem 1.1 is then equivalent to

$$\mathbb{P}(A) \leq c\rho^k + \mathbb{P}(B)$$

for constants $0 < \rho < 1$ and $c$ that depend only on $s, t$.

The proof proceeds as follows. Fix a large $N$ and let

$$A_i^N = \{(S_{i,i+n})_{n=1}^N \text{ has } k \text{ upcrossings of } (s,t)\}.$$

Since $A = \bigcup A_1^N$, it suffices to show that, for each $N$,

$$\mathbb{P}(A_1^N) \leq c\rho^k + \mathbb{P}(B)$$

with $c, \rho$ independent of $N$. Fix an integer $R$ much bigger than $N$ (we will eventually take $R \to \infty$) and let $I \subseteq \{1, \ldots, R\}$ be the random set of indices defined by

$$I = \{i \in [1; R] : A_i^N \text{ occurs}\}.$$

By stationarity of $S_{i,j}$, we have

$$\mathbb{P}(A_1^N) = \frac{1}{R} \sum_{i=1}^R \mathbb{P}(A_i^N) = \frac{1}{R} \mathbb{E}\left( \sum_{i=1}^R \chi_{A_i^N} \right) = \frac{1}{R} \mathbb{E}|I|.$$

We proceed to estimate the expected size of $I$. We divide $I$ into two parts:

$$I_0 = \{i \in I : B_i \text{ occurs}\} \quad \text{and} \quad I_1 = I \setminus I_0.$$

Since $\frac{1}{R}\mathbb{E}(|I_0|) = \mathbb{P}(B)$, it suffices to show that $\frac{1}{R}|I_1| \leq c\rho^k$.

By definition, for each $i \in I_1$, there is a (random) sequence of $k$ pairs of nonempty intervals $U_i(1) \subseteq V_i(1) \subseteq \cdots \subseteq U_i(k) \subseteq V_i(k)$ whose left endpoint is $i$ and with length $\leq N$, and such that $S_{U_i(m)} < s$ and $S_{V_i(m)} > t$ for $1 \leq m \leq k$.

We now pass to a subsequence of the $U_i$'s and $V_i$'s by performing two refinements of the sequence. First, clearly, $|U_i(m+1)| > |V_i(m)| > |U_i(m)|$ because $S_{U_{i+1}(m)}, S_{U_i(m)} \neq S_{V_i(m)}$, and because of the given inclusions. Thus, $|U_i(m)| > 2(m-1)$. If we delete the first $k_0 = \lceil k/2 \rceil + 1$ pairs, then we are left with a sequence $U_i(k_0+1) \subseteq V_i(k_0+1) \subseteq \cdots \subseteq U_i(k) \subseteq V_i(k)$ of at least $k' = \lfloor k/2 \rfloor - 1$ pairs of intervals, all of which are of length greater than $k$.

Second, notice that if $|V_i(m)|/|U_i(m)| \leq 1 + \delta$, then $i \in I_0$ since $U_i(m)$ would $\delta$-fill $V_i(m)$. Thus, for $i \in I_1$, we also have $|V_i(m)|/|U_i(m)| > 1 + \delta$. Choose $q$ so that

$$(1+\delta)^{q-1} \geq 72/\delta^2,$$

that is, $q = \lceil \log(72/\delta^2)/\log(1+\delta) \rceil + 1$. By deleting the intervals $U_i(j), V_i(j+1)$ when $j \neq 0 \pmod{q}$ and renumbering the remaining ones, we are left with a sequence of $k'' \geq \lfloor (k'-1)/q \rfloor$ pairs of intervals $\widetilde{U}_i(1) \subseteq \widetilde{V}_i(1) \subseteq \cdots \subseteq \widetilde{U}_i(k'') \subseteq \widetilde{V}_i(k'')$, all having length $> k$ and satisfying $S_{\widetilde{U}_i(m)} < s$, $S_{\widetilde{V}_i(m)} > t$, and which, additionally, satisfy the growth condition $|V_i(m)| \geq \frac{72}{\delta^2}|U_i(m)|$ for every $1 \leq m \leq k''$.

We now apply a combinatorial result whose proof we defer to Section 6.

LEMMA 5.1. *Let $\varepsilon < 1/4$. Suppose that $J \subseteq \mathbb{N}$ is finite and that for each $j \in J$, we are given a sequence of intervals $\widetilde{U}_j(1) \subseteq \widetilde{V}_j(1) \subseteq \cdots \subseteq \widetilde{U}_j(L) \subseteq \widetilde{V}_j(L)$ with left endpoint $j$ and satisfying $|\widetilde{V}_j(n)| \geq \frac{2}{\varepsilon^2}|\widetilde{U}_j(n)|$. Suppose that none of the $\widetilde{V}_j(n)$'s can be $6\varepsilon$-filled by a disjoint collection of $\widetilde{U}_i(m)$'s. Then,*

$$|J| \leq \left(1 + \frac{\varepsilon}{6}\right)^{-(L-1)/(\log 1/\varepsilon)} \cdot \left|\bigcup_{j \in J} \widetilde{V}_j(L)\right|.$$

We apply the lemma to our situation with $J = I_1$, $L = k''$ and $\varepsilon = \delta/6$. The hypothesis is satisfied by definition of $I_1$. It follows that there are constants $c > 0$ and $0 < \rho < 1$ depending only on $\delta$ (hence on $s, t$), such that

$$|I_1| \leq c \cdot \rho^{-k} \cdot \left|\bigcup_{i \in I_1} V_i(k'')\right|$$

$$\leq c \cdot \rho^{-k}(N+R)$$

and the last inequality holds because $\bigcup_{i \in I_1} V_i(k'') \subseteq [1; N+R]$. We thus have

$$\mathbb{P}(A) = \frac{1}{R}\mathbb{E}|I| = \frac{1}{R}\mathbb{E}|I_0| + \frac{1}{R}\mathbb{E}|I_1| \leq \mathbb{P}(B) + c \cdot \rho^{-k}\left(1 + \frac{N}{R}\right)$$

and the proof of Theorem 1.1 is complete by taking $R \to \infty$.

**6. Proof of Lemma 5.1.** The remainder of this paper is devoted to the proof of Lemma 5.1. Some of the statements below are standard; we supply proofs for completeness. Others parts of the argument are related to the effective Vitali covering lemma from [10]. See [12, 14] for other examples of covering lemmas in probability and ergodic theory.

We say that a collection of segments is *disjoint* if its members are pairwise disjoint. The following is a version of the classical Vitali covering lemma:

LEMMA 6.1. *If $\mathcal{V}$ is a collection of intervals, then there is a disjoint subcollection $\mathcal{V}' \subseteq \mathcal{V}$ with $|\bigcup \mathcal{V}'| \geq |\bigcup \mathcal{V}|/2$.*

PROOF. Let $\mathcal{V}' \subseteq \mathcal{V}$ be a minimal collection satisfying $\bigcup \mathcal{V}' = \bigcup \mathcal{V}$. Order the intervals in $\mathcal{V}'$ by their left endpoints, say $\mathcal{V}' = \{V_1, V_2, \ldots, V_m\}$. Then, the subsequence consisting of intervals with even indices is disjoint, similarly for the subsequence with odd indices. One of these must cover at least half of $\mathcal{V}$. □

For $\varepsilon > 0$ and an interval $U = [a; b]$, the $\varepsilon$-blowup of $U$ is

$$U^\varepsilon = [a - \varepsilon|U|, b + \varepsilon|U|] \cap \mathbb{Z}.$$

Note that $U \subseteq U^\varepsilon$ and $|U^\varepsilon| \leq (1 + 2\varepsilon)|U|$. For a collection $\mathcal{U}$ of intervals, we write $\mathcal{U}^\varepsilon = \{U^\varepsilon : U \in \mathcal{U}\}$.

LEMMA 6.2. *If $\mathcal{U}$ is a collection of intervals, then $|\bigcup \mathcal{U}^\varepsilon| \leq (1+2\varepsilon)|\bigcup \mathcal{U}|$.*

PROOF. Let $A = \bigcup \mathcal{U}$ and decompose $A$ into disjoint maximal intervals $V_1, \ldots, V_k$, so $|\bigcup \mathcal{U}| = \sum |V_i|$. For each $V_i$, one clearly has

$$\bigcup_{U \in \mathcal{U}: U \subseteq V_i} U^\varepsilon \subseteq V_i^\varepsilon,$$

thus

$$\left|\bigcup \mathcal{U}^\varepsilon\right| \leq \sum_{i=1}^{k} \left|\bigcup_{U \in \mathcal{U}: U \subseteq V_i} U^\varepsilon\right| \leq \sum_{i=1}^{k} (1+2\varepsilon)|V_i| = (1+2\varepsilon)\left|\bigcup \mathcal{U}\right|. \quad □$$

A *tower* of height $M$ over a finite set $I \subseteq \mathbb{N}$ is a collection $\mathcal{U} = \{U_i(k) : i \in I, 1 \leq k \leq M\}$ of intervals such that $i$ is the left endpoint of $U_i(k)$ [we shall actually only use the fact that $i \in U_i(k)$] and for each $i \in I$, the sequence $U_i(1) \subseteq U_i(2) \subseteq \cdots \subseteq U_i(M)$ is strictly increasing. The $k$th *level* of $\mathcal{U}$ is the collection

$$\mathcal{U}(k) = \{U_i(k) : i \in I\}.$$

Note that the intervals in $\mathcal{U}(k)$ are not necessarily of the same size and, although $|U_i(k)| < |U_i(k+1)|$, it need not be true that $|U_i(k)| \leq |U_j(k+1)|$ if $i \neq j$.

Let $\mathcal{U} = \{U_i(k)\}$ be a tower of height $M$ over a set $I$. The $\varepsilon$-crust of $\mathcal{U}$ is the set of $V \in \mathcal{U}$ whose $\varepsilon$-blowup is strictly maximal with respect to inclusion, that is,

$$\mathcal{V} = \{V \in \mathcal{U}(M) : \text{if } V^\varepsilon \subsetneq W^\varepsilon \text{ for some } W \in \mathcal{U}, \text{ then } V = W\}.$$

It is clear that $\bigcup \mathcal{U} \subseteq \bigcup \mathcal{V}^\varepsilon$.

LEMMA 6.3. *Let $0 < \varepsilon < 1$. Suppose that $\mathcal{U} = \{U_i(k)\}$ is a tower over $I$ of height $2$ satisfying*

$$|U_i(2)| \geq \frac{2}{\varepsilon^2}|U_i(1)| \qquad \text{for all } i \in I$$

*and $\mathcal{V} \subseteq \mathcal{U}(2)$ is the $\varepsilon$-crust of $\mathcal{U}$. Then:*

(1) *for each $U \in \mathcal{U}(1)$ and $V \in \mathcal{V}$, if $U \cap V \neq \varnothing$ then $U \subseteq V^\varepsilon$;*
(2) *there exists $\widehat{\mathcal{U}} \subseteq \mathcal{U}(1)$ and a disjoint $\widehat{\mathcal{V}} \subseteq \mathcal{V}$ such that $|\bigcup \widehat{\mathcal{U}}| \leq \frac{1}{2}|\bigcup \mathcal{U}|$, $(\bigcup \widehat{\mathcal{U}}) \cap (\bigcup \widehat{\mathcal{V}}) = \varnothing$ and $\bigcup \mathcal{U} \subseteq (\bigcup \widehat{\mathcal{V}}^\varepsilon) \cup (\bigcup \widehat{\mathcal{U}})$.*

PROOF. Let $U \in \mathcal{U}(1)$ and $V \in \mathcal{V}$ with $U \cap V \neq \varnothing$. In order to show that $U \subseteq V^\varepsilon$, it suffices to show that $\varepsilon|V| \geq |U|$. Let $i \in I$ be such that $U = U_i(1)$ and write $W = U_i(2)$, the interval "above" $U$ in $\mathcal{U}$. Since $V$ is in the $\varepsilon$-crust, we cannot have $V^\varepsilon \subsetneq W^\varepsilon$. Since $V \cap W \neq \varnothing$, this implies that

$$(1+\varepsilon)|V| \geq |W| \geq \frac{2}{\varepsilon^2}|U|,$$

which gives the desired conclusion.

To establish (2), apply the Vitali lemma to $\mathcal{V}$ to obtain a disjoint family $\widehat{\mathcal{V}} \subseteq \mathcal{V}$ with $|\bigcup \widehat{\mathcal{V}}| \geq \frac{1}{2}|\bigcup \mathcal{V}|$. Let

$$\widehat{\mathcal{U}} = \left\{U \in \mathcal{U}(1) : U \cap \left(\bigcup \widehat{\mathcal{V}}\right) = \varnothing\right\}.$$

The conclusion now follows from (1). $\square$

We want to replace the constant $1/2$ in the Vitali lemma with a constant close to 1. This can be achieved using the standard trick of applying the Vitali lemma to several layers of covers and iteratively disjointifying each level in turn.

Henceforth, all logarithms are taken to base 2.

LEMMA 6.4. *Let $0 < \varepsilon < 1$. Suppose that $\mathcal{U} = \{U_j(k)\}$ is a tower of height $M \geq 1 + \log(1/\varepsilon)$ over a set $J$ and*

$$|U_j(k+1)| > \frac{2}{\varepsilon^2}|U_j(k)|.$$

*There is then a disjoint subcollection $\mathcal{V} \subseteq \mathcal{U}$ such that $|\bigcup \mathcal{V}| \geq (1-3\varepsilon)|\bigcup \mathcal{U}|$.*

PROOF. Set $\mathcal{U}_0 = \mathcal{U}$ and $\mathcal{V}_0 = \varnothing$. For $1 \leq n < M$, we inductively define subcollections $\mathcal{V}_n, \mathcal{U}_n \subseteq \mathcal{U}(M-n)$ satisfying:

(1) $\mathcal{V}_n, \mathcal{U}_n \subseteq \mathcal{U}_{n-1}$;
(2) $\bigcup \mathcal{U}_{n-1} \subseteq (\bigcup \mathcal{U}_n) \cup (\bigcup \mathcal{V}_n^\varepsilon)$;
(3) $(\bigcup \mathcal{V}_n) \cap (\bigcup \mathcal{U}_n) = \varnothing$;
(4) $|\bigcup \mathcal{U}_n| \leq \frac{1}{2}|\bigcup \mathcal{U}_{n-1}|$.

To produce $\mathcal{U}_n, \mathcal{V}_n$, we apply part (2) of Lemma 6.3 to the top two layers of $\mathcal{U}_{n-1}$. Clearly, the collection $\mathcal{V} = \bigcup_{1 \leq k < M} \mathcal{V}_k$ is disjoint and $\bigcup \mathcal{U} \subseteq (\bigcup \mathcal{V}^\varepsilon) \cup (\bigcup \mathcal{U}_{M-1})$. By property (4), we have

$$\left|\bigcup \mathcal{U}_{M-1}\right| \leq \left(\frac{1}{2}\right)^{M-1}\left|\bigcup \mathcal{U}\right| \leq \varepsilon \left|\bigcup \mathcal{U}\right|,$$

so

$$\left|\bigcup \mathcal{V}^\varepsilon\right| \geq \left|\bigcup \mathcal{U}\right| - \left|\bigcup \mathcal{U}_{M-1}\right| \geq (1-\varepsilon)\left|\bigcup \mathcal{U}\right|.$$

Using the inequality $|\bigcup \mathcal{V}| \geq \frac{1}{1+2\varepsilon}|\bigcup \mathcal{V}^\varepsilon|$ from Lemma 6.2, we have

$$\left|\bigcup \mathcal{V}\right| \geq \frac{1-\varepsilon}{1+2\varepsilon}\left|\bigcup \mathcal{U}\right| \geq (1-3\varepsilon)\left|\bigcup \mathcal{U}\right|,$$

as desired. □

For the remainder of this section, we adopt the following notation. Fix an integer $L$ and $0 < \varepsilon < 1$, a finite set $I \subseteq \mathbb{Z}$ and two towers $\mathcal{U} = \{U_i(k)\}$ and $\mathcal{V} = \{V_i(k)\}$ of height $L+1$ over $I$ satisfying

$$U_i(0) \subseteq V_i(0) \subseteq U_i(1) \subseteq V_i(1) \subseteq \cdots \subseteq U_i(L) \subseteq V_i(L)$$

(for convenience, we start from level 0) and

$$|V_i(k)| \geq \frac{2}{\varepsilon^2}|U_i(k)|.$$

Note that this ensures a similar growth rate for the substack $\mathcal{U}$.

LEMMA 6.5. *Let $L \geq 1 + \log(1/\varepsilon)$. Then, either there is a $V \in \mathcal{V}(L)$ which can be $6\varepsilon$-filled by a disjoint subcollection of $\mathcal{U}$, or else $|\bigcup \mathcal{V}(L)| \geq (1+\frac{\varepsilon}{6})|\bigcup \mathcal{U}(0)|$.*

PROOF. Select a maximal disjoint subset $\mathcal{W}$ of the $\varepsilon$-crust of $\mathcal{V}$. By the Vitali lemma and Lemma 6.3, we have

$$\sum_{W \in \mathcal{W}} |W| = \left|\bigcup \mathcal{W}\right| \geq \frac{1}{2(1+2\varepsilon)} \left|\bigcup \mathcal{V}(L)\right| \geq \frac{1}{2(1+2\varepsilon)} \left|\bigcup \mathcal{U}(0)\right|.$$

We distinguish between two cases. First, if every $W \in \mathcal{W}$ satisfies $|W \setminus \bigcup \mathcal{U}(L-1)| > \varepsilon |W|$, then, since $\mathcal{W}$ is disjoint, we would have

$$\left|\bigcup \mathcal{V}\right| \geq \left|\bigcup \mathcal{U}(0)\right| + \left|\bigcup \mathcal{W} \setminus \bigcup \mathcal{U}(L-1)\right|$$
$$= \left|\bigcup \mathcal{U}(0)\right| + \sum_{W \in \mathcal{W}} \left|W \setminus \bigcup \mathcal{U}(L-1)\right|$$
$$\geq \left|\bigcup \mathcal{U}(0)\right| + \sum_{W \in \mathcal{W}} \varepsilon |W|$$
$$\geq \left(1 + \frac{\varepsilon}{2(1+2\varepsilon)}\right) \left|\bigcup \mathcal{U}(0)\right|,$$

which gives the desired bound.

Otherwise, let $W \in \mathcal{W}$ be such that $|W \setminus \bigcup \mathcal{U}(L-1)| < \varepsilon |W|$. Let

$$\mathcal{Y} = \{U_i(k) : 0 \leq k \leq L-1 \text{ and } U_i(L-1) \cap W \neq \varnothing\},$$
$$\mathcal{Z} = \{U_i(k) : 0 \leq k \leq L-1 \text{ and } U_i(L-1) \subseteq W\}.$$

By assumption, we know that

(6.1) $$\left|W \cap \left(\bigcup \mathcal{Y}\right)\right| \geq (1-\varepsilon)|W|.$$

We claim that then

(6.2) $$\left|W \cap \left(\bigcup \mathcal{Z}\right)\right| \geq (1-3\varepsilon)|W|.$$

Suppose that this were not the case. There is then a subset $A \subseteq W$ of size $> \varepsilon |W|$, not covered by $\mathcal{Z}$ and all of whose points are at a distance of at least $\varepsilon |W|$ from $\mathbb{Z} \setminus W$. By inequality (6.1), there is some $U_i(L-1) \in \mathcal{Y} \setminus \mathcal{Z}$ that intersects $A$ at a point $r$. Since $U_i(L-1) \notin \mathcal{Z}$, it must intersect $\mathbb{Z} \setminus W$ at a point $r'$. Since $U_i(L-1)$ contains the interval with endpoints $r, r'$, we have

$$|U_i(L-1)| \geq |r - r'| \geq \varepsilon |W|,$$

so

$$|V_i(L-1)| \geq \frac{2}{\varepsilon^2} |U_i(L-1)| > 2|W|.$$

On the other hand, $U_i(L-1) \cap W \neq \varnothing$ [because $U_i(L-1) \in \mathcal{Y}$] and since $U_i(L-1) \subseteq V_i(L-1)$, we also have $V_i(L-1) \cap W \neq \varnothing$. Since $W$ is in the $\varepsilon$-crust of $\mathcal{V}$, Lemma 6.3 implies that $V_i(L-1) \subseteq W^\varepsilon$, contradicting the size bound we got for $V_i(L-1)$.

To complete the proof, we apply Lemma 6.4 to the tower $\mathcal{Z}$. We obtain a disjoint subcollection of $\mathcal{Z}$ (and hence of $\mathcal{U}$) whose members are contained in $W$ and have total size at least $(1-3\varepsilon)|\cup \mathcal{Z}|$, which, by inequality (6.2), is at least $(1-6\varepsilon)|W|$, as required. $\square$

We can now prove Lemma 5.1, which we rephrase as follows (notice that our tower is now numbered starting at 0 and that we have removed the tildes from the notation).

LEMMA 6.6. *With the above notation, suppose that no interval $V \in \mathcal{V}$ can be $6\varepsilon$-filled by disjoint elements of $\mathcal{U}$. Then,*

$$\left|\bigcup \mathcal{U}(0)\right| \leq \left(1+\frac{\varepsilon}{6}\right)^{-\lfloor L/\log(1/\varepsilon) \rfloor} \left|\bigcup \mathcal{V}(L)\right|.$$

PROOF. Set $M = \log(1/\varepsilon)$. It suffices to prove

$$\left|\bigcup \mathcal{U}\right| \geq \left(1+\frac{\varepsilon}{6}\right)^{-\lfloor L/M \rfloor} \left|\bigcup \mathcal{U}(0)\right|$$

and it is enough to prove this when $L$ is an integer multiple of $M$; write $L = kM$. We proceed by induction on $k$. The base of the induction is the previous lemma. Now, given that it is true for $k$ and given $L = (k+1)M$, we can apply the induction hypothesis to the restrictions of $\mathcal{U}, \mathcal{V}$ to levels $0, 1, \ldots, kM$. This tells us that

$$\left|\bigcup \mathcal{U}(kM)\right| \geq \left(1+\frac{\varepsilon}{6}\right)^k \left|\bigcup \mathcal{U}(0)\right|.$$

Now, consider the restriction of the towers to levels $kM, kM+1, \ldots, (k+1)M$. Applying the base case, we get

$$\left|\bigcup \mathcal{U}((k+1)M)\right| \geq \left(1+\frac{\varepsilon}{6}\right) \left|\bigcup \mathcal{U}(kM)\right|.$$

Putting these together completes the proof. $\square$

As mentioned in the Introduction, everything above can be carried out for symmetric intervals and for cubes in $\mathbb{Z}^d$; the proofs generalize easily to that case, although the constants change. We note that the Vitali lemma (Lemma 6.1) requires a different proof in higher dimensions, but this is classical.

For completeness, we provide the proof of the higher-dimensional analog of Lemma 6.2. Consider the case of squares in $\mathbb{Z}^2$. The $\varepsilon$-blowup of a square $U \times V$ is $U^\varepsilon \times V^\varepsilon$, which can be written as a disjoint union,

$$U^\varepsilon \times V^\varepsilon = (U \times U) \cup B_1 \cup B_2 \cup B_3 \cup B_4,$$

where $B_1 = U^\varepsilon \times V \setminus U \times V$ are two vertical strips of width $\varepsilon$, $B_2 = U \times V^\varepsilon \setminus U \times V$ are two horizontal strips of height $\varepsilon$, $B_3$ is the union of two $\varepsilon \times \varepsilon$ squares outside the upper-left and lower-right corners of $U \times V$ and $B_4$ is the union of two $\varepsilon \times \varepsilon$ squares outside the upper-right and lower-left corners of $U \times V$. To obtain an analog of Lemma 6.2, we must show that if $\{U_i \times V_i\}$ is a collection of squares and $U_i^\varepsilon \times V_i^\varepsilon = (U_i \times V_i) \cup \bigcup_{t=1,2,3,4} B_{i,t}$ as above, then, for each $t = 1, 2, 3, 4$,

$$\left|\bigcup_i B_{i,t}\right| \le 2\varepsilon \left|\bigcup_i U_i \times V_i\right|.$$

This follows from the one-dimensional case by decomposing $\bigcup_i (U_i \times V_i \cup B_{i,t})$ into the union of the intersection of this set with parallel translates of lines. For instance, for $t = 1$, the intersection of $\bigcup_i (U_i \times V_i \cup B_{i,t})$ with each horizontal line is the $\varepsilon$-blowup (in the one-dimensional sense) of the intersection of $\bigcup_i U_i \times V_i$ with that line and so the one-dimensional lemma can be applied. Now, sum over all lines.

The proof for cubes in $\mathbb{Z}^d$ is proved by induction on the dimension using a similar strategy.

**Acknowledgments.** This work was done as part of the author's Ph.D. studies under the guidance of Benjamin Weiss, whom I would like to thank for all his support and advice. I also thank the anonymous referee for a careful reading and for suggesting many simplifications in the proofs.

EINSTEIN INSTITUTE OF MATHEMATICS
EDMOND J. SAFRA CAMPUS
GIVAT RAM
THE HEBREW UNIVERSITY OF JERUSALEM
JERUSALEM 91904
ISRAEL
E-MAIL: hochman@math.princeton.edu